\documentclass[12pt,a4paper]{amsart} 
\usepackage{amssymb,amsfonts,amsmath,amsthm,amscd}
\usepackage[latin1]{inputenc}

\newtheorem{theorem}{Theorem}[section]
\newtheorem*{theorem/maintheorem}{Theorem \ref{main-theorem}}

\newtheorem{lemma}[theorem]{Lemma}
\newtheorem{corollary}[theorem]{Corollary}
\newtheorem{example}[theorem]{Example}
\newtheorem{remark}[theorem]{Remark}
\newtheorem{discussion}[theorem]{Discussion}

\numberwithin{equation}{section}

\newcommand{\fm}{\mbox{$\mathfrak{m}$}}

\newcommand{\cm}{\mbox{$\mathcal{M}$}}

\newcommand{\rees}{\mbox{$\mathbf{R}$}}
\newcommand{\graded}{\mbox{$\mathbf{G}$}}
\newcommand{\fiber}{\mbox{$\mathbf{F}$}}
\newcommand{\symmetric}{\mbox{$\mathbf{S}$}}

\newcommand{\reltype}{\mbox{$\rm{rt}$}}
\newcommand{\rednumber}{\mbox{$\rm{r}$}}
\newcommand{\spread}{\mbox{$\textit{l}$}}
\newcommand{\height}{\mbox{$\rm{height}$}}
\newcommand{\grade}{\mbox{$\rm{grade}$}}

\newcommand{\lcm}{\mbox{$\rm{lcm}$}}

\newcommand{\koszul}[3]{\mbox{$H_{#1}(#2\,;#3)$}}
\newcommand{\cycle}[3]{\mbox{$Z_{#1}(#2\,;#3)$}}
\newcommand{\boundary}[3]{\mbox{$B_{#1}(#2\,;#3)$}}

\begin{document}

\title{The equations of Rees algebras of equimultiple ideals of deviation one}

\author{Ferran Mui\~nos and Francesc Planas-Vilanova}
\thanks{The second author is partially supported by the MICINN Grant MTM2010-20279.}

\begin{abstract}
We describe the equations of the Rees algebra $\rees(I)$ of an
equimultiple ideal $I$ of deviation one provided that $I$ has a
reduction generated by a regular sequence $x_1,\ldots ,x_s$ such that
the initial forms $x_1^*,\ldots ,x_{s-1}^*$ are a regular sequence in
the associated graded ring. In particular, we prove that there is a
single equation of maximum degree in a minimal generating set of the
equations of $\rees(I)$, which recovers some previous known results.
\end{abstract}


\maketitle

\section{Introduction}\label{introduction}

Let $(R,\fm)$ be a Noetherian local ring and let $I$ be an
equimultiple ideal of $R$ of deviation $\mu(I)-\height(I)=1$. The aim
of this paper is to study the equations of
$\rees(I)=R[It]=\oplus_{n\geq 0}I^n$, the Rees algebra of $I$. We will
assume that $x_1,\ldots,x_s,y$ is a minimal generating set of $I$,
that $J=(x_1,\ldots,x_s)$ is a reduction of $I$ generated by an
$R$-sequence and that the initial forms $x_1^*,\ldots,x_{s-1}^*$ in
$\graded(I)=\oplus_{n\geq 0}I^{n}/I^{n+1}$, the associated graded ring
of $I$, are a $\graded(I)$-sequence.

Recall that the reduction number of $I$ with respect to $J$ is the
least integer $r\geq 0$ such that $I^{r+1}=JI^r$, denoted by
$\rednumber_J(I)$ (notice that, with our assumptions, $r\geq 1$). Set
$V=R[X_1,\ldots,X_s,Y]$ and let $\varphi:V\to\rees(I)$ be the
polynomial presentation of $\rees(I)$ sending $X_i$ to $x_it$ and $Y$
to $yt$. Let $Q=\oplus_{n\geq 1}Q_n$ be the kernel of $\varphi$, whose
elements will be referred to as the equations of $\rees(I)$. For
instance, for $r=\rednumber_{J}(I)$, the containment of $y^{r+1}$ in
$JI^{r}$ induces an equation $Y^{r+1}-\sum X_iF_i\in Q_{r+1}$, with
$F_i\in V_r$. Given an integer $n\geq 1$, set $Q\langle
n\rangle\subset Q$ the ideal generated by the homogeneous equations of
$\rees(I)$ of degree at most $n$ in $X_1,\ldots,X_s,Y$.  The relation
type of $I$, denoted by $\reltype(I)$, is the least integer $N\geq 1$
such that $Q=Q\langle N\rangle$. When $Q=Q\langle 1\rangle$, $I$ is
called of linear type. Our main result is the following.

\begin{theorem/maintheorem}
Let $(R,\fm)$ be a Noetherian local ring and let $I$ be an ideal of
$R$. Let $x_1,\ldots,x_s,y$ be a minimal generating set of $I$ and
$J=(x_1,\ldots,x_s)$. Assume that $x_1,\ldots,x_s$ is an $R$-sequence
and that $x_1^*,\ldots,x_{s-1}^*$ is a $\graded(I)$-sequence. Then,
for each $n\geq 2$, the map sending $F\in Q_{n}$ to
$F(0,\ldots,0,1)\in (JI^{n-1}:y^{n})$ induces an isomorphism of
$R$-modules
\begin{eqnarray*}
\left[ \frac{Q}{Q\langle n-1\rangle}\right]_{n}\cong
\frac{JI^{n-1}:y^{n}}{JI^{n-2}:y^{n-1}}.
\end{eqnarray*}
In particular, if $J$ is a reduction of $I$ with reduction number
$r=\rednumber_J(I)$, then $\reltype(I)=\rednumber_J(I)+1$ and there is
a form $Y^{r+1}-\sum X_iF_i\in Q_{r+1}$, with $F_i\in V_r$, such that
$Q=(Y^{r+1}-\sum X_iF_i)+Q\langle r\rangle$.
\end{theorem/maintheorem}

Note that, with these assumptions, $s\leq \grade(I)\leq
\height(I)\leq\spread(I)\leq\mu(J)=s$, where $\spread(\cdot)$ and $\mu(\cdot)$
stand for analytic spread and minimal number of generators,
respectively. Thus $I$ is an equimultiple ideal (i.e.,
$\height(I)=l(I)$) of deviation $\mu(I)-\height(I)=1$. Roughly
speaking, the theorem says how to obtain a minimal generating set of
the equations of $\rees(I)$. For the equations of degree 1: pick a
minimal generating set of the first syzygies of $I$, viewed as
elements of $Q_{1}$; for the equations of higher degree $n$, $2\leq
n\leq r+1$, take representatives of the inverse images of a minimal
generating set of $(JI^{n-1}:y^{n})/(JI^{n-2}:y^{n-1})$ (see
Remark~\ref{min-gen-set-of-equations} and
Example~\ref{exemple-classic}).

Thus far, the study of the equations of $\rees(I)$ has produced a vast
literature. Part of this work has been focused on ideals having small
deviations as well as on the interplay between the reduction number
and the relation type (see, just to mention a few of them,
\cite{huckaba1}, \cite{huckaba2}, \cite{hh1}, \cite{hh2},
\cite{trung1}, \cite{trung2}, \cite{vasconcelos1},
\cite{vasconcelos2}). The particular hypotheses and interests in this
note owe much to the works of T.~Cortadellas and S.~Zarzuela in
\cite{cz}, W.~Heinzer and M.-K.~Kim in \cite{hk}, S.~Huckaba in
\cite{huckaba1} and \cite{huckaba2}, N.V.~Trung in \cite{trung1} and
\cite{trung2}, and W.V.~Vasconcelos in \cite{vasconcelos1}. In fact,
Theorem \ref{main-theorem} sprouted as an attempt to understand
\cite[Theorem 2.3.3]{vasconcelos1}. That is one of the main reasons
for considering ideals $I$ of the form $I=(J,y)$, $J$ being a
reduction of $I$ generated by a sufficiently good sequence. The reader
may also consult \cite{hsv1}, \cite{hsv2} for a recent account on the
equations of $\rees(I)$ with similar assumptions.

Even in a simpler case, i.e. $\mu(I)-\height(I)=0$, the equations of
$\rees(I)$ may be difficult to describe. Remarkably, any ideal $I$ of
the principal class is generated by an $R$-sequence $x_{1},\ldots
,x_{s}$, provided that $I$ is prime (\cite{davis}) or $R$ is
Cohen-Macaulay. In both cases the equations are generated by the
Koszul relations $x_iX_j-x_jX_i$, $1\leq i<j\leq s$. In particular,
$I$ is of linear type. However, if $I$ is not prime and $R$ is not
Cohen-Macaulay, this is no longer true. For instance, consider the
ideals generated by a system of parameters. C.~Huneke asked in
\cite{huneke} whether there is a uniform bound for the relation type
of these ideals in a complete local equidimensional Noetherian ring
$R$. The full answer to this question was given in \cite{wang},
\cite{lai} and \cite{agh}.  Concretely, in \cite[Example 2.1]{agh}, it
was shown that if the non-Cohen-Macaulay locus of $R$ has dimension 2
or more, there exist families of parametric ideals of $R$ with
unbounded relation type. This gives an idea of the complexity of the
structure of the equations of $\rees(I)$.

In our case, i.e. $\mu(I)-\height(I)=1$, if $R$ is Cohen-Macaulay, one
has that $I$ is of linear type if and only if $I$ is locally of the
principal class at all minimal primes of $I$ (see
\cite[Theorem~4.8]{hmv}). But, if $I$ is not locally of the principal
class at its minimal primes, then the relation type may be arbitrarily
large. For instance, take $R=k[[a_{1},a_{2}]]$, the power series ring
in two variables $a_1,a_2$ over a field $k$, set $x_1=a_1^p$,
$x_2=a_2^p$ and $y=a_1a_2^{p-1}$, with $p\geq 2$. Let
$I=(x_{1},x_{2},y)$.  Then $I$ has deviation one and is
$(a_{1},a_{2})$-primary, hence not locally of the principal class at
its unique minimal prime. One can check that $\reltype(I)\geq p$. In
fact, $I$ fulfills the hypotheses of Theorem~\ref{main-theorem},
$J=(x_{1},x_{2})$ being a reduction of $I$. Thus, the containment
$y^p\in JI^{p-1}$, induces the only equation of degree $p$ in a
minimal generating set of equations of $\rees(I)$ and $\reltype(I)=p$
(see Example~\ref{exemple-classic}).

As said before, this note has its origins in the following result of
W.V.~Vasconcelos in \cite[Theorem~2.3.3]{vasconcelos1}: let $(R,\fm)$
be a Cohen-Macaulay local ring of dimension $d$ and let $I$ be an
$\fm$-primary ideal of $R$. Let $x_{1},\ldots, x_{d},y$ be a minimal
generating set of $I$, where $J=(x_{1},\ldots ,x_{d})$ is a reduction
of $R$ with reduction number $\rednumber_J(I)=1$. Then there is a form
$Y^{2}-\sum X_iF_i\in Q_{2}$, with $F_i\in V_1$, such that
$Q=(Y^{2}-\sum X_iF_i)+Q\langle 1\rangle$. In particular,
$\reltype(I)=2$.

Remark that the hypotheses of \cite[Theorem~2.3.3]{vasconcelos1} imply
that $x_{1},\ldots, x_{d}$ is an $R$-sequence and that the initial
forms $x_1^*,\ldots,x_d^*$ are a $\graded(I)$-sequence (see the result
of P.~Valabrega and G.~Valla in \cite[Proposition~3.1]{vv}). By
Theorem~\ref{main-theorem}, it is enough to suppose
$x_1^*,\ldots,x_{d-1}^*$ being a $\graded(I)$-sequence and to consider
any reduction number (see Corollary~\ref{crelle-generalisation}).

Also as a corollary of Theorem~\ref{main-theorem}, we recover the
result of W.~Heinzer and M.-K.~Kim in \cite[Theorem~5.6]{hk}, where
they prove that the equations of $\fiber(I)=\rees(I)\otimes
R/\fm=\oplus_{n\geq 0}I^n/\fm I^n$, the fiber cone of $I$, are
generated by a unique equation of degree $\rednumber_J(I)+1$.

The paper is organised as follows: in Section~\ref{equations-of-u} we
recall from \cite{planas} how to express the equations of a standard
algebra in terms of the effective relations and the Koszul homology.
Section~\ref{equations-of-rees(I)} is focused on the equations of the
Rees algebra. In Section~\ref{main-result} we prove
Theorem~\ref{main-theorem} and the aforementioned results of
W.V.~Vasconcelos and W.~Heinzer and M.-K.~Kim. We finish by giving
some examples in Section~\ref{some-examples}.

All the unexplained notations can be found in \cite{bh} and \cite{sh}.
All along the paper, $(R,\fm)$ is a Noetherian local ring and
$J\subset I$ are two ideals of $R$.

\section{The equations of a standard algebra and the effective
relations}\label{equations-of-u}

By a standard $R$-algebra $U$ we mean a graded $R$-algebra
$U=\oplus_{n\geq 0}U_{n}$, with $U_{0}=R$, $U=R[U_{1}]$ and $U_{1}$
minimally generated by $z_{1},\ldots ,z_{s}\in U_1$ as an
$R$-module. For instance, the $R$-algebra $\rees(I)=\oplus_{n\geq
  0}I^{n}$, the $R/I$-algebra $\graded(I)=\oplus_{n\geq
  0}I^{n}/I^{n+1}$ and the $R/\fm$-algebra $\fiber(I)=\oplus_{n\geq
  0}I^n/\fm I^n$ are standard algebras.

Let $V=R[T_1,\ldots ,T_s]$ be a polynomial ring with variables
$T_{1},\ldots ,T_s$ and let $\varphi:V\rightarrow U$ be the induced
presentation of $U$ sending $T_{i}$ to $z_{i}$. Let $Q=\oplus_{n\geq
  1}Q_{n}$ be the kernel of $\varphi$, whose elements will be referred
to as the equations of $U$. As before, set $Q\langle n\rangle$ the
ideal generated by the homogeneous equations of $U$ of degree at most
$n$. The relation type of $U$, denoted by $\reltype(U)$, is the least
integer $N\geq 1$, such that $Q=Q\langle N\rangle$. Observe that the
relation type of $I$, $\reltype(I)$, defined in
Section~\ref{introduction}, is nothing else but the relation type of
its Rees algebra $\rees(I)$.

Although $Q$ depends on the presentation $\varphi$, the quotients
$(Q/Q\langle n-1\rangle)_{n}=Q_{n}/V_{1}Q_{n-1}$, $n\geq 2$, do not.
Indeed, let $\symmetric(U_{1})$ be the symmetric algebra of $U_{1}$
and let $\alpha:\symmetric(U_{1})\rightarrow U$ be the canonical
morphism induced by the identity in degree one. Given $n\geq 2$, the
module of effective $n$-relations of $U$ is defined to be
$E(U)_{n}=\ker(\alpha_{n})/U_{1}\ker(\alpha_{n-1})$. One can prove
that, for any $n\geq 2$, $(Q/Q\langle n-1\rangle)_{n}\cong E(U)_{n}$
(see \cite[Definition~2.2]{planas}). In particular, the relation type
of $U$ can be calculated as the least integer $N\geq 1$, such that
$E(U)_{n}=0$ for all $n\geq N+1$.

This description of $(Q/Q\langle n-1\rangle)_{n}$ as
$E(U)_{n}=\ker(\alpha_{n})/U_{1}\ker(\alpha_{n-1})$ has the advantage
of being canonical. However, it is often useful to express
$(Q/Q\langle n-1\rangle)_{n}$ as a Koszul homology module (see
\cite[Corollary~2.7]{planas}). Concretely, $(Q/Q\langle
n-1\rangle)_{n}\cong E(U)_{n}\cong \koszul{1}{z}{U}_{n}$, for $n\geq
2$, where $\koszul{1}{z}{U}_n$ denotes the first homology module of
the complex:
\begin{eqnarray*}
\ldots \to \wedge_{2}(R^{s})\otimes
U_{n-2}\stackrel{\partial_{2,n-2}}{\longrightarrow}
\wedge_{1}(R^{s})\otimes
U_{n-1}\stackrel{\partial_{1,n-1}}{\longrightarrow}U_{n}\to 0,
\end{eqnarray*}
where the Koszul differentials are defined as follows: if $e_1,\ldots
,e_s$ stands for the canonical basis of $R^s$ and $u\in U_{n-2}$ and
$v\in U_{n-1}$, then
\begin{eqnarray*}
\partial_{2,n-2}(e_i\wedge e_j\otimes u)=e_j\otimes z_i\cdot
u-e_i\otimes z_j\cdot u \; \mbox{ and } \; \partial_{1,n-1}(e_i\otimes
v)=z_i\cdot v.
\end{eqnarray*}
As usual, $\cycle{1}{z}{U}$ and $\boundary{1}{z}{U}$ will stand for
the graded modules of $1$-cycles and $1$-boundaries, respectively, of
the Koszul complex of $z$. Observe that $\boundary{1}{z}{U}_{1}=0$.

For the sake of easy reference we finish this section by stating the
following two remarks.

\begin{remark}\label{explicit-presentation} {\rm
For $n\geq 2$, there is an isomorphism $\tau_n: (Q/Q\langle
n-1\rangle)_{n}\stackrel{\cong}{\to}\koszul{1}{z}{U}_{n}$ sending the
class of $F\in Q_n$ modulo $Q\langle n-1 \rangle$ to the homology
class of the $s$-tuple $(F_1(z),\ldots,F_s(z))\in
\cycle{1}{z}{U}_{n}\subset U_{n-1}\oplus \ldots \oplus U_{n-1}$, where
$F_1,\ldots,F_s$ are elements in $V_{n-1}$ satisfying
$F=T_1F_1+\ldots+T_sF_s$. For $n=1$, setting $Q\langle 0\rangle =0$
and taking into account that $\boundary{1}{z}{U}_{1}=0$, then there is
an isomorphism $\tau_1:Q_1\to\cycle{1}{z}{U}_{1}$.  }\end{remark}

\begin{proof}
Take $f=\varphi:V\to U$ and $g={\rm id}_{V}$ in
\cite[Theorem~2.4]{planas}.
\end{proof}

\begin{remark}\label{min-gen-set-of-equations} {\rm 
Let $U$ be a standard $R$-algebra and let $Q=\oplus_{n\geq 1}Q_{n}$ be
the equations of $U$.  Suppose that $Q=Q\langle N\rangle$ for some
$N\geq 1$. Set
\begin{eqnarray*}
\mathcal{C}=\{ F_{1,1},\ldots ,F_{1,s_{1}},\ldots ,F_{N,1},\ldots
,F_{N,s_{N}}\}, 
\end{eqnarray*}
with $F_{i,j}\in Q_{i}$. Then $\mathcal{C}$ is a minimal generating
set of $Q$ if and only if, for each $n=1,\ldots ,N$, the classes of
$F_{n,1},\ldots ,F_{n,s_{n}}$ modulo $Q\langle n-1\rangle$ are a
minimal generating set of $(Q/Q\langle n-1\rangle)_{n}$.
}\end{remark}

\section{The equations of Rees algebras via Koszul homology}
\label{equations-of-rees(I)} 

As always, $(R,\fm)$ is a Noetherian local ring and $I$ is an ideal of
$R$. The purpose of this section is to describe the modules
$(Q/Q\langle n-1\rangle)_n$, where $Q$ are the equations of
$\rees(I)$. By Remark~\ref{explicit-presentation}, this is equivalent
to describe the graded components of a Koszul homology module.
Although the next lemma is written for $\rees(I)$, similar results can
also be established for $\graded(I)$ or $\fiber(I)$.

\begin{lemma}\label{succexcurta}
Let $(R,\fm)$ be a Noetherian local ring and let $I$ be an ideal of
$R$. Let $x_1,\ldots,x_s,y$ be a minimal generating set of $I$ and
$J=(x_1,\ldots,x_s)$. Then, for each integer $n\geq 2$, there is a
short exact sequence
\begin{eqnarray*}
0\to \frac{\koszul{1}{x_1t,\ldots,x_st}{\rees(I)}_n}
{yt\koszul{1}{x_1t,\ldots,x_st}{\rees(I)}_{n-1}} & \longrightarrow
&\\ \longrightarrow \koszul{1}{x_1t,\ldots,x_st,yt}{\rees(I)}_n &
\stackrel{\sigma_n}{\longrightarrow} &
\frac{JI^{n-1}:y^n}{JI^{n-2}:y^{n-1}} \to 0.
\end{eqnarray*}
Moreover, $\sigma_n$ sends the class of a cycle
$(w_1t^{n-1},\ldots,w_st^{n-1},w_{s+1}t^{n-1})$, $w_{i}\in I^{n-1}$,
to the class of $a\in (JI^{n-1}:y^{n})$, where $w_{s+1}=ay^{n-1}+b$
for some $b\in JI^{n-2}$.
\end{lemma}

\begin{proof}
Take $z=x_1t,\ldots,x_st$ and $\underline{z}=x_1t,\ldots,x_st,yt$ in
$It\subset \rees(I)=R[It]=\oplus_{n\geq 0}I^n$. Consider the induced
graded long exact sequence of Koszul homology:
\begin{eqnarray*}
\koszul{1}{z}{\rees(I)}_{n-1}\stackrel{(\rho_1)_n}{\longrightarrow}
\koszul{1}{z}{\rees(I)}_n & \longrightarrow & \\ \longrightarrow
\koszul{1}{\underline{z}}{\rees(I)}_n\longrightarrow
\koszul{0}{z}{\rees(I)}_{n-1} & \stackrel{(\rho_0)_n}{\longrightarrow}
& \koszul{0}{z}{\rees(I)}_n,
\end{eqnarray*}
where $(\rho_i)_n$ is just the multiplication by $\pm yt$. We get the
following short exact sequence:
\begin{eqnarray*}
0 \to\textrm{coker}(\rho_1)_n \longrightarrow
\koszul{1}{\underline{z}}{\rees(I)}_n
\stackrel{\sigma^{\prime}_n}{\longrightarrow}\ker(\rho_0)_n \to 0.
\end{eqnarray*}
Clearly, ${\rm coker}(\rho_1)_n=
\koszul{1}{z}{\rees(I)}_n/yt\koszul{1}{z}{\rees(I)}_{n-1}$. On the
other hand,
\begin{eqnarray*}
\koszul{0}{z}{\rees(I)}_{n-1}=I^{n-1}/JI^{n-2}.
\end{eqnarray*}
 Thus $\ker(\rho_0)_n=I^{n-1}\cap (JI^{n-1}:y)/JI^{n-2}$. One can
 check that $\sigma^{\prime}_n$ maps the homology class of a cycle
 $(w_1t^{n-1},\ldots,w_st^{n-1},w_{s+1}t^{n-1})\in
 \oplus_{i=1}^{s+1}I^{n-1}t^{n-1}$ to the class of $w_{s+1}\in
 I^{n-1}\cap (JI^{n-1}:y)$ modulo $JI^{n-2}$. Finally, consider the
 mapping $\theta_n$,
\begin{eqnarray*}
\ker(\rho_0)_n=\frac{I^{n-1}\cap(JI^{n-1}:y)}{JI^{n-2}}
\stackrel{\theta_n}{\longrightarrow}\frac{JI^{n-1}:y^n}{JI^{n-2}:y^{n-1}},
\end{eqnarray*}
defined as follows: for each $w\in I^{n-1}\cap(JI^{n-1}:y)$, since
$I^{n-1}=y^{n-1}R+JI^{n-2}$, take $a\in R$ and $b\in JI^{n-2}$ such
that $w=ay^{n-1}+b$. Clearly, $a\in (JI^{n-1}:y^{n})$. Let
$\overline{w}$ be the class of $w$ modulo $JI^{n-2}$ and let
$\overline{a}$ the class of $a$ modulo $(JI^{n-2}:y^{n-1})$. Set
$\theta_n(\overline{w})=\overline{a}$. An elementary computation shows
that $\theta_n$ is a well-defined isomorphism. Take $\sigma_n=
\theta_n\circ\sigma^{\prime}_n$.
\end{proof}

Although Lemma~\ref{succexcurta} has a version for $n=1$, we are more
interested in a kind of analogue obtained from considering
$x_{1},\ldots ,x_{s}$ and $x_1,\ldots ,x_s,y$ as sequences of elements
of degree zero in $I\subset R$.

\begin{remark}\label{syzygies} {\rm 
Let $(R,\fm)$ be a Noetherian local ring and let $I$ be an ideal of
$R$. Let $x_1,\ldots,x_s,y$ be a minimal generating set of $I$ and
$J=(x_1,\ldots,x_s)$. Then there is an exact sequence:
\begin{eqnarray*}
0\to\frac{\koszul{1}{x_1,\ldots,x_s}{R}}
{y\koszul{1}{x_1,\ldots,x_s}{R}}\longrightarrow
\koszul{1}{x_1,\ldots,x_s,y}{R}\stackrel{\tilde{\sigma}_1}{\longrightarrow}
\frac{(J:y)}{J}\to 0,
\end{eqnarray*}
where $\tilde{\sigma}_1$ sends the homology class of a cycle
$(w_1,\ldots,w_s,w_{s+1})\in R^{s+1}$ to the class of $w_{s+1}\in
(J:y)$. In particular, if $x_1,\ldots,x_s$ is an $R$-sequence,
$\koszul{1}{x_1,\ldots,x_s}{R}=0$ and $\tilde{\sigma}_1$ is an
isomorphism. }\end{remark}

In the next lemma we will prove that, for $n\geq 2$,
$\koszul{1}{x_1t,\ldots,x_st}{\rees(I)}_{n}=0$, provided that
$x_{1},\ldots ,x_{s}$ is an $R$-sequence such that $x_{1}^*,\ldots
,x_{s-1}^*$ is a $\graded(I)$-sequence. Notice that we are not saying
that the whole $\koszul{1}{x_1t,\ldots,x_st}{\rees(I)}$ vanishes. In
fact, if $s>1$, $x_{1}(x_{2}t)\in x_{1}t\cdot\rees(I)$ whereas
$x_{1}\not\in x_{1}t\cdot\rees(I)$, hence $x_1t,\ldots,x_st$ is not an
$\rees(I)$-sequence.

Before this, recall the result of P.~Valabrega and G.~Valla in
\cite[Corollary~2.7]{vv} which characterizes to be a
$\graded(I)$-sequence.  Let $J\subset I$ be two ideals of a Noetherian
local ring $R$. Let $x_1,\ldots ,x_s$ be a minimal generating set of
$J$ and write $J_{0}=0$ and $J_{i}=(x_{1},\ldots ,x_{i})$, for
$i=1,\ldots ,s$. In particular, the initial forms $x_1^*,\ldots
,x_s^*$ in $\graded(I)$ are in $I/I^2$. Then, $x_1^*,\ldots ,x_s^*$ is
a $\graded(I)$-sequence if and only if $x_1,\ldots ,x_s$ is an
$R$-sequence and the so-called Valabrega-Valla modules
$VV_{J_i}(I)_{n}=J_i\cap I^n/J_iI^{n-1}$ are zero for all $i=1,\ldots
,s$ and all $n\geq 1$.

\begin{lemma}\label{homologia-zero}
Let $(R,\fm)$ be a Noetherian local ring and let $J\subset I$ be two
ideals of $R$. Let $x_1,\ldots,x_s$ be a minimal generating set of
$J$. Write $J_{0}=0$ and $J_{i}=(x_{1},\ldots ,x_{i})$, for
$i=1,\ldots ,s$. Let $n\geq 2$. Then the following two conditions are
equivalent:
\begin{itemize}
\item[$(i)$] $\koszul{1}{x_1t,\ldots,x_it}{\rees(I)}_n=0$, for all
  $i=1,\ldots ,s$;
\item[$(ii)$] $(J_{i-1}I^{n-1}:x_i)\cap I^{n-1}/J_{i-1}I^{n-2}=0$, for
  all $i=1,\ldots ,s$.
\end{itemize}
Suppose that, in addition, $x_1,\ldots,x_s$ is an $R$-sequence and
$x_1^*,\ldots,x_{s-1}^*$ is a $\graded(I)$-sequence. Then, for all
$n\geq 2$ and all $i=1,\ldots ,s$,
\begin{eqnarray*}
\koszul{1}{x_1t,\ldots,x_it}{\rees(I)}_n=0.
\end{eqnarray*}
\end{lemma}
\begin{proof}
For a fixed $n\geq 2$, let us prove the equivalence between $(i)$ and
$(ii)$ by induction on $s\geq 1$. If $s=1$, then
\begin{eqnarray*}
\koszul{1}{x_1t}{\rees(I)}_{n}\cong (0:x_1)\cap I^{n-1}\cong
(J_0I^{n-1}:x_1)\cap I^{n-1}/J_0I^{n-2},
\end{eqnarray*}
and the claim follows.

Now, take $s>1$ and set $z=x_1t,\ldots,x_{s-1}t$ and
$\underline{z}=x_1t,\ldots,x_{s-1}t,x_{s}t$. The induced graded long
exact sequence of Koszul homology gives rise to the exact sequence:
\begin{eqnarray*}
\koszul{1}{z}{\rees(I)}_n\longrightarrow\koszul{1}{\underline{z}}{\rees(I)}_n
\longrightarrow\koszul{0}{z}{\rees(I)}_{n-1}\stackrel{(\rho_0)_n}{\longrightarrow}
\koszul{0}{z}{\rees(I)}_n,
\end{eqnarray*}
where $(\rho_0)_n$ is the multiplication by $\pm x_st$. Then
\begin{eqnarray*}
\ker(\rho_{0})_{n}=
\frac{((x_1t,\ldots,x_{s-1}t):x_{s}t)_{n-1}}{(x_1t,\ldots,x_{s-1}t)_{n-1}}
\cong \frac{(J_{s-1}I^{n-1}:x_{s})\cap I^{n-1}}{J_{s-1}I^{n-2}}.
\end{eqnarray*}
Therefore, one has the following exact sequence:
\begin{eqnarray*}
\koszul{1}{z}{\rees(I)}_n\longrightarrow\koszul{1}{\underline{z}}{\rees(I)}_n
\longrightarrow\frac{(J_{s-1}I^{n-1}:x_{s})\cap I^{n-1}}{J_{s-1}I^{n-2}}\to 0.
\end{eqnarray*}
By the induction hypothesis, $\koszul{1}{x_1t,\ldots
  ,x_it}{\rees(I)}_n=0$, for all $i=1,\ldots ,s-1$, is equivalent to
$(J_{i-1}I^{n-1}:x_i)\cap I^{n-1}/J_{i-1}I^{n-2}=0$, for all
$i=1,\ldots ,s-1$. Therefore, using the former exact sequence, one
deduces that $\koszul{1}{x_1t,\ldots ,x_it}{\rees(I)}_n=0$, for all
for all $i=1,\ldots ,s$, is equivalent to $(J_{i-1}I^{n-1}:x_i)\cap
I^{n-1}/J_{i-1}I^{n-2}=0$, for all $i=1,\ldots ,s$. 

Suppose that $x_1,\ldots,x_s$ is an $R$-sequence and
$x_1^*,\ldots,x_{s-1}^*$ is a $\graded(I)$-sequence. Let us see that
$(J_{i-1}I^{n-1}:x_{i})\cap I^{n-1}/J_{i-1}I^{n-2}=0$, for all
$i=1,\ldots ,s$ and all $n\geq 2$. Using the aforementioned result of
P.~Valabrega and G.~Valla,
\begin{eqnarray*}
VV_{J_i}(I)_{n-1}=J_i\cap I^{n-1}/J_iI^{n-2}=0
\end{eqnarray*}
for all $i=1,\ldots ,s-1$ and all $n\geq 2$.  Since $x_1,\ldots,x_s$
is an $R$-sequence, then $(J_{i-1}:x_{i})=J_{i-1}$ for all
$i=1,\ldots,s$. Therefore, for all $i=1,\ldots, s$ and all $n\geq 2$,
\begin{eqnarray*}
J_{i-1}I^{n-2}\subseteq (J_{i-1}I^{n-1}:x_{i})\cap I^{n-1}\subseteq
(J_{i-1}:x_{i})\cap I^{n-1}=J_{i-1}\cap I^{n-1}=J_{i-1}I^{n-2}.
\end{eqnarray*}
\end{proof}

One can state a different version of the second part of
Lemma~\ref{homologia-zero}, which will in turn lead to a slightly
different version of the main result of this paper (see
Remark~\ref{different-version-theorem}).

Recall that $x_{1},\ldots ,x_{s}$ is a $d$-sequence if
$(J_i:x_{i+1}x_{j})=(J_i:x_{j})$ for all $0\leq i\leq s-1$ and all
$j\geq i+1$ (where $J_0=0$, $J_{i}=(x_1,\ldots ,x_i)$ and $J_s=J$, as
before). This condition is equivalent to $(J_i:x_{i+1})\cap J= J_i$
for all $0\leq i\leq s-1$. Clearly, $R$-sequences are $d$-sequences.

\begin{remark}\label{different-version-lemma} {\rm
Let $(R,\fm)$ be a Noetherian local ring and let $I$ be an ideal of
$R$. Let $x_1,\ldots,x_s$ be a minimal generating set of $J$, where
$J=(x_1,\ldots,x_s)$ is a reduction of $I$ with reduction number
$r=\rednumber_J(I)$. Assume that
\begin{itemize}
\item[$(a)$] $x_1,\ldots,x_s$ is a $d$-sequence and
\item[$(b)$] $VV_{J_i}(I)_{r+1}=(x_{1},\ldots ,x_{i})\cap
  I^{r+1}/(x_{1},\ldots ,x_{i})I^{r}=0$ for all $i=1,\ldots ,s-1$.
\end{itemize}
Then $\koszul{1}{x_1t,\ldots,x_it}{\rees(I)}_n=0$, for all $n\geq r+2$
and all $i=1,\ldots ,s$. Suppose that, in addition,
\begin{itemize}
\item[$(c)$] $x_1,\ldots,x_s$ is an $R$-sequence and
\item[$(d)$] $VV_{J_{i}}(I)_{r}=(x_{1},\ldots ,x_{i})\cap
  I^{r}/(x_{1},\ldots ,x_{i})I^{r-1}=0$ for all $i=1,\ldots ,s-1$.
\end{itemize}
Then $\koszul{1}{x_1t,\ldots,x_it}{\rees(I)}_{r+1}=0$, for all
$i=1,\ldots ,s$.  }\end{remark}
\begin{proof}
Since $J$ is a reduction of $I$ with reduction number $r$,
$I^{r+1}=JI^{r}$. Using $(a)$ and $(b)$, for $i=1,\ldots ,s$, we get:
\begin{eqnarray*}
(J_{i-1}:x_{i})\cap I^{r+1}=(J_{i-1}:x_i)\cap J\cap
  I^{r+1}=J_{i-1}\cap I^{r+1}=J_{i-1}I^{r}.
\end{eqnarray*}
Therefore, for all $i=1,\ldots ,s$, $(J_{i-1}:x_{i})\cap
I^{r+1}=J_{i-1}I^{r}$.  Now, using the result of N.V.~Trung in
\cite[Proposition~4.7(i)]{trung2}, one has $(J_{i-1}:x_{i})\cap
I^{n}=J_{i-1}I^{n-1}$, for all $n\geq r+1$ and all $i=1,\ldots
,s$. This clearly implies that
\begin{eqnarray*}
(J_{i-1}I^{n}:x_{i})\cap I^{n}=J_{i-1}I^{n-1},
\end{eqnarray*}
for all $n\geq r+1$ and all $i=1,\ldots ,s$, which by
Lemma~\ref{homologia-zero}, is equivalent to
$$\koszul{1}{x_1t,\ldots,x_it}{\rees(I)}_n=0,$$ for all $n\geq r+2$ and
all $i=1,\ldots ,s$. 

Suppose now $(c)$ and $(d)$. Then, for all $i=1,\ldots ,s$,
\begin{eqnarray*}
J_{i-1}I^{r-1}\subseteq (J_{i-1}I^{r}:x_i)\cap I^{r}\subseteq
(J_{i-1}:x_i)\cap I^r=J_{i-1}\cap I^r=J_{i-1}I^{r-1}.
\end{eqnarray*}
Therefore, for all $i=1,\ldots ,s$, $(J_{i-1}I^{r}:x_i)\cap
I^{r}=J_{i-1}I^{r-1}$. By Lemma~\ref{homologia-zero},
$$\koszul{1}{x_1t,\ldots,x_it}{\rees(I)}_{r+1}=0,$$ for all $i=1,\ldots
,s$.
\end{proof}

\section{Main result}\label{main-result}

We have now all the ingredients to prove the main result of the
paper. As always, set $V=R[X_1,\ldots,X_s,Y]$ and $Q$ the kernel of
the polynomial presentation $\varphi:V\to\rees(I)$ sending $X_i$ to
$x_it$ and $Y$ to $yt$.

\begin{theorem}\label{main-theorem}
Let $(R,\fm)$ be a Noetherian local ring and let $I$ be an ideal of
$R$. Let $x_1,\ldots,x_s,y$ be a minimal generating set of $I$ and
$J=(x_1,\ldots,x_s)$. Assume that $x_1,\ldots,x_s$ is an $R$-sequence
and that $x_1^*,\ldots,x_{s-1}^*$ is a $\graded(I)$-sequence. Then,
for each $n\geq 2$, the map sending $F\in Q_{n}$ to
$F(0,\ldots,0,1)\in (JI^{n-1}:y^{n})$ induces an isomorphism of
$R$-modules
\begin{eqnarray*}
\left[ \frac{Q}{Q\langle n-1\rangle}\right]_{n}\cong
\frac{JI^{n-1}:y^{n}}{JI^{n-2}:y^{n-1}}.
\end{eqnarray*}
In particular, if $J$ is a reduction of $I$ with reduction number
$r=\rednumber_J(I)$, then $\reltype(I)=\rednumber_J(I)+1$ and there is
a form $Y^{r+1}-\sum X_iF_i\in Q_{r+1}$, with $F_i\in V_r$, such that
$Q=(Y^{r+1}-\sum X_iF_i)+Q\langle r\rangle$.
\end{theorem}
\begin{proof}
By Lemma~\ref{succexcurta}, with $z=x_1t,\ldots,x_st$ and
$\underline{z}=x_1t,\ldots,x_st,yt$,
\begin{eqnarray*}
0\to\frac{\koszul{1}{z}
  {\rees(I)}_{n}}{yt\koszul{1}{z}{\rees(I)}_{n-1}}\longrightarrow
\koszul{1}{\underline{z}}{\rees(I)}_n\stackrel{\sigma_n}{\longrightarrow}
\frac{JI^{n-1}:y^{n}}{JI^{n-2}:y^{n-1}} \to 0
\end{eqnarray*}
is an exact sequence for all $n\geq 2$. By Lemma~\ref{homologia-zero},
$\koszul{1}{z}{\rees(I)}_{n}=0$ for all $n\geq 2$. Therefore
$\koszul{1}{\underline{z}}{\rees(I)}_n\stackrel{\sigma_n}{\cong}
(JI^{n-1}:y^{n})/(JI^{n-2}:y^{n-1})$, for $n\geq 2$. With
Remark~\ref{explicit-presentation}, we conclude:
\begin{eqnarray*}
\left[ \frac{Q}{Q\langle n-1\rangle}\right]_{n}
\stackrel{\tau_n}{\cong}\koszul{1}{\underline{z}}{\rees(I)}_n
\stackrel{\sigma_n}{\cong} \frac{JI^{n-1}:y^{n}}{JI^{n-2}:y^{n-1}},
\end{eqnarray*}
for all $n\geq 2$. Given $F\in Q_{n}$, write
$F=\sum_{i=1}^{s}X_{i}F_{i}+YG$, with $F_{i},G\in V_{n-1}$. Then the
morphism $\tau_n$ sends the class of $F$ to the homology class of
$(F_{1}(\underline{z}),\ldots ,F_{s}(\underline{z}),
G(\underline{z}))$. But $G(\underline{z})=G(x_{1},\ldots
,x_{s},y)t^{n-1}$ and $G(x_{1},\ldots ,x_{s},y)=G(0,\ldots
,0,1)y^{n-1}+b$, for some $b\in JI^{n-2}$. By Lemma~\ref{succexcurta},
$\sigma_n$ sends the homology class of $(F_{1}(\underline{z}),\ldots
,F_{s}(\underline{z}), G(\underline{z}))$ to the class of
$G(0,\ldots,0,1)$ modulo $(JI^{n-2}:y^{n-1})$ and notice that
$G(0,\ldots,0,1)=F(0,\ldots,0,1)$.

If $J$ is a reduction of $I$ with reduction number $r$,
$(JI^{n-1}:y^{n})=R$ for all $n\geq r+1$. Therefore $(Q/Q\langle
n-1\rangle)_n=0$ for all $n>r+1$ and $(Q/Q\langle r\rangle)_{r+1}\cong
R/(JI^{r-1}:y^r)\neq 0$, since $y^r\not\in JI^{r-1}$. Therefore
$\reltype(I)=r+1$. Finally, notice that the containment $y^{r+1}\in
JI^{r}$ induces an equation of the form $Y^{r+1}-\sum_iX_iF_i$, with
$F_i\in V_r$, which is sent by $\sigma_n\circ\tau_n$ to the class of
$1$ in $R/(JI^{r-1}:y^r)$.
\end{proof}

\begin{remark}\label{weaker-hyp} {\rm 
In Theorem~\ref{main-theorem}, we use the hypotheses
``$x_1,\ldots,x_s$ is an $R$-sequence and $x_1^*,\ldots,x_{s-1}^*$ is
a $\graded(I)$-sequence'' to assure
$\koszul{1}{x_1t,\ldots,x_st}{\rees(I)}_{n}=0$, for all $n\geq 2$, by
means of Lemma~\ref{homologia-zero}. Therefore, in
Theorem~\ref{main-theorem}, one can substitute, if needed, the
hypotheses ``$x_1,\ldots,x_s$ is an $R$-sequence and
$x_1^*,\ldots,x_{s-1}^*$ is a $\graded(I)$-sequence'' by the weaker
hypotheses ``$(J_{i-1}I^{n-1}:x_i)\cap I^{n-1}/J_{i-1}I^{n-2}=0$, for
all $i=1,\ldots ,s$ and all $n\geq 2$''. See Example~\ref{gen-classic}
as an application of this comment.  }\end{remark}

With the weaker hypotheses pinpointed in
Remark~\ref{different-version-lemma}, we get the next version of
Theorem~\ref{main-theorem}.

\begin{remark}\label{different-version-theorem} {\rm
Let $(R,\fm)$ be a Noetherian local ring and let $I$ be an ideal of
$R$. Let $x_1,\ldots,x_s,y$ be a minimal generating set of $I$, where
$J=(x_1,\ldots,x_s)$ is a reduction of $I$ with reduction number
$r=\rednumber_J(I)$. Assume that
\begin{itemize}
\item[$(a)$] $x_1,\ldots,x_s$ is a $d$-sequence and
\item[$(b)$] $VV_{J_i}(I)_{r+1}=(x_{1},\ldots ,x_{i})\cap
  I^{r+1}/(x_{1},\ldots ,x_{i})I^{r}=0$ for all $i=1,\ldots ,s-1$.
\end{itemize}
Then $\reltype(I)\leq \rednumber_J(I)+1$. Suppose that, in addition,
\begin{itemize}
\item[$(c)$] $x_1,\ldots,x_s$ is an $R$-sequence and
\item[$(d)$] $VV_{J_{i}}(I)_{r}=(x_{1},\ldots ,x_{i})\cap
  I^{r}/(x_{1},\ldots ,x_{i})I^{r-1}=0$ for all $i=1,\ldots ,s-1$.
\end{itemize}
Then $\reltype(I)=\rednumber_J(I)+1$ and there is a form $Y^{r+1}-\sum
X_iF_i\in Q_{r+1}$, with $F_i\in V_r$, such that $Q=(Y^{r+1}-\sum
X_iF_i)+Q\langle r\rangle$.  }\end{remark}
\begin{proof}
It follows from the proof of Theorem~\ref{main-theorem}, but using
Remark~\ref{different-version-lemma} instead of
Lemma~\ref{homologia-zero}.
\end{proof}

\begin{discussion}\label{general-hypotheses}{\rm 
The hypotheses of Remark~\ref{different-version-theorem} are connected
with the works of S.~Huckaba in \cite[Theorem~1.4]{huckaba2} and
N.V.~Trung in \cite[Theorem~6.4]{trung2} (see also
\cite[Theorem~3.2]{cz}, \cite{huckaba1}, \cite[Theorem~5.3]{gp},
\cite{trung1}). In \cite[Theorems~1.4, 1.5]{huckaba2}, S.~Huckaba
proved that if $I$ is an ideal with $\spread(I)=\height(I)+1\geq 2$
and such that any minimal reduction $J$ of $I$ can be generated by a
$d$-sequence $x_{1},\ldots ,x_{s}$ with $x_{1}^{*},\ldots
,x_{s-1}^{*}$ being a $\graded(I)$-sequence ($s=\spread(I)$), then
$\reltype(I)\leq \rednumber_J(I)+1$. If in addition
$\mu(I)=\spread(I)+1$, then the equality $\reltype(I)=
\rednumber_J(I)+1$ holds.  In particular, $r=\rednumber_J(I)$ is
independent of $J$. In fact, N.V.~Trung improved this last result in
\cite[Theorem~6.4]{trung2} by showing that $r$ coincides with the
Castelnuovo-Mumford regularity of $\rees(I)$. To prove
$\reltype(I)\geq\rednumber_J(I)+1$, S.~Huckaba showed that the
equality $I^{r+1}=JI^r$ induces an equation of $\rees(I)$ of maximum
degree. In Theorem~\ref{main-theorem} and with a different approach,
we have completed the description of the whole ideal of equations of
$\rees(I)$.}\end{discussion}

As a corollary of Theorem \ref{main-theorem}, one can prove the result
of W.V.~Vasconcelos in \cite[Theorem~2.3.3]{vasconcelos1} for any
reduction number not necessarily equal to $1$.

\begin{corollary}\label{crelle-generalisation}
Let $(R,\fm)$ be a Cohen-Macaulay local ring of dimension $d$. Let $I$
be an $\fm$-primary ideal of $R$ minimally generated by $x_{1},\ldots,
x_{d},y$, where $J=(x_{1},\ldots
,x_{d})$ is a reduction of $R$ with reduction number
$r=\rednumber_J(I)$. Suppose that $x_{1}^{*},\ldots ,x_{d-1}^{*}$ is a
$\graded(I)$-sequence. Then there is a form $Y^{r+1}-\sum X_iF_i\in
Q_{r+1}$, with $F_i\in V_r$, such that $Q=(Y^{r+1}-\sum
X_iF_i)+Q\langle r\rangle$.  In particular, $\reltype(I)=r+1$.
\end{corollary}
\begin{proof}
Since $R$ is Cohen-Macaulay and $I$ is $\fm$-primary, 
$x_{1},\ldots ,x_{d}$ is an $R$-sequence and the results follows from
Theorem~\ref{main-theorem}.
\end{proof}

Before stating the next result, recall that for an ideal $L$ of $R$
and any standard $R$-algebra $U$, the relation type of $U\otimes R/L$
is $\reltype(U\otimes R/L)\leq \reltype(U)$ (see
e.g. \cite[Example~3.2]{planas}). Hence $\reltype(\fiber(I))\leq
\reltype(\graded(I))\leq \reltype(I)$. In fact, for any $n\geq 2$ and
writing $E(I)_n$ for $E(\rees(I))_n$ (see Section~\ref{equations-of-u}
to recall the definition of effective relations), there is an exact
sequence $E(I)_{n+1}\to E(I)_{n}\to E(\graded(I))_{n}\to 0$. In
particular, $\reltype(I)=\reltype(\graded(I))$ and $E(I)_{N}\cong
E(\graded(I))_{N}$ for $N=\reltype(I)$
(\cite[Proposition~3.3]{planas}; see also \cite[p.~268]{hku}). In our
context, $\reltype(\fiber(I))$ is also equal to $\reltype(I)$, which
is not a general fact (see the Conjecture of G.~Valla in
\cite[\S~2]{hmv} and a counterexample in \cite[Example~4.4]{suv}; see
also \cite[p.~268 and Corollary~2.6]{hku}).

As a corollary of Theorem~\ref{main-theorem}, we recover the
following result of W.~Heinzer and M.-K.~Kim in
\cite[Theorem~5.6]{hk}.

\begin{corollary}\label{hk-generalisation}
Let $(R,\fm)$ be a Noetherian local ring with infinite residue field
$k=R/\fm$ and let $I$ be an ideal of $R$. Let $x_1,\ldots,x_s,y$ be a
minimal generating set of $I$, where $J=(x_1,\ldots,x_s)$ is a
reduction of $I$ with reduction number $r=\rednumber_J(I)$. Assume
that $x_1,\ldots,x_s$ is an $R$-sequence and that
$x_1^*,\ldots,x_{s-1}^*$ is a $\graded(I)$-sequence. Then there is a
form $Y^{r+1}-\sum X_iF_i$, with $F_i\in k[X_1,\ldots ,X_s,Y]$ forms
of degree $r$ and $\fiber(I)\cong k[X_1,\ldots ,X_s,Y]/(Y^{r+1}-\sum
X_iF_i)$. In particular,
$\reltype(\fiber(I))=\rednumber_J(I)+1=\reltype(I)$.
\end{corollary}

\begin{proof}
By Theorem~\ref{main-theorem}, $\reltype(\fiber(I))\leq
\reltype(I)=r+1$. By \cite[Lemma~5.2]{hk}, $E(\fiber(I))_{n}=0$ for
all $2\leq n\leq r$ and $E(\fiber(I))_{r+1}\neq 0$. Thus $\fiber(I)$
has only equations of degree $r+1$ and $\reltype(\fiber(I))=r+1$. By
Theorem~\ref{main-theorem}, $E(I)_{r+1}$ is cyclic and generated by
the equation of $\rees(I)$ induced by the containment $y^{r+1}\in
JI^{r}$. Therefore the same happens with $E(\fiber(I))_{r+1}$ (see
\cite[Proposition~3.2]{gp} or \cite[p.~268]{hku} ).
\end{proof}

\section{Some examples}\label{some-examples}

Our purpose is to take advantage of Theorem~\ref{main-theorem} to
obtain a minimal generating set of the equations of $\rees(I)$. The
ascending chain of colon ideals $\{(JI^{n-1}:y^{n})\}_{n\geq 1}$,
which need not be rigid but stabilizes in $R$, can be calculated in
any computer algebra system, giving a possibly alternative procedure
to find the equations of $\rees(I)$.

The following ideal $I$, for a specific $p\geq 1$, is often used as an
example of an ideal of relation type at least $p$. However, we are not
aware of any reference with a detailed description of a minimal
generating set of the equations of $\rees(I)$, for a general $p\geq
1$. As said above, Theorem~\ref{main-theorem} will be crucial to our
purposes.

\begin{example}\label{exemple-classic} {\rm
Let $(R,\fm)$ be a Noetherian local ring. Let $a_1,a_2$ an
$R$-sequence and $p\geq 2$. Set $x_{1}=a_{1}^{p}$, $x_{2}=a_{2}^{p}$
and $y=a_{1}a_{2}^{p-1}$. Let $I$ be the ideal generated by
$x_{1},x_{2},y$. Set $V=R[X_{1},X_{2},Y]$ and let $\varphi:V\to
\rees(I)$ be the presentation of $\rees(I)$ sending $X_{i}$ to $x_it$
and $Y$ to $yt$. Then a minimal generating set of the ideal
$Q=\ker(\varphi)$ is obtained from:
\begin{itemize}
\item a unique equation
  $F_n(X_1,X_2,Y)=a_1^{p-n}Y^n-a_2^{p-n}X_1X_2^{n-1}\in Q_{n}$ of
  degree $n$, for each $n$, $2\leq n\leq p$;
\item two equations $F_1(X_1,X_2;Y)=a_1^{p-1}Y-a_2^{p-1}X_1$ and
  $G_1(X_1,X_2,Y)=a_2Y-a_1X_2\in Q_{1}$ of degree 1.
\end{itemize}
}\end{example}
\begin{proof}
\noindent We start by proving the hypotheses in
Theorem~\ref{main-theorem}. Clearly $x_{1},x_{2}$ is an $R$-sequence
and $J=(x_1,x_2)$ is a reduction of $I$ since $I^{p}=JI^{p-1}$. By
\cite[Corollary~3]{ks}, a monomial $m$ on $a_{1},a_{2}$ belongs to an
ideal generated by monomials $m_{1},\ldots ,m_{r}$ on $a_{1},a_{2}$ if
and only if $m$ is a multiple of some $m_{i}$. It follows that
$y^{p-1}\not\in JI^{p-2}$ and $I^{p-1}\not\subseteq JI^{p-2}$. Thus
$\rednumber_J(I)=p-1$.

\noindent \underline{\sc Claim}: $x_{1}^{*}$ is a
$\graded(I)$-sequence.  By \cite[Corollary~2.7]{vv}, it suffices to
prove 
\begin{eqnarray*}
VV_{(x_1)}(I)_n=x_{1}R\cap I^{n}/x_{1}I^{n-1}=0
\end{eqnarray*} 
for all $n\geq 1$. Fix $n\geq 1$. By \cite[Proposition~1]{ks},
\begin{eqnarray*}
x_{1}R\cap I^{n}=(L_{i,j,k}\; |\; 
i,j,k \text{ positive integers such that } i+j+k=n),
\end{eqnarray*} 
where
$L_{i,j,k}=\lcm(a_{1}^{p},a_{1}^{ip}a_{2}^{jp}(a_{1}a_{2}^{p-1})^{k})$.
Let us prove that $L_{i,j,k}$ is in $x_1I^{n-1}$.

Indeed, if $i\geq 1$, then
$L_{i,j,k}=a_{1}^{ip}a_{2}^{jp}(a_{1}a_{2}^{p-1})^{k}=
a_{1}^{p}[a_{1}^{(i-1)p}a_{2}^{jp}(a_{1}a_{2}^{p-1})^{k}]\in
x_{1}I^{n-1}$ and we have finished. Hence we can suppose $i=0$ and
$j+k=n$. If $k=0$, then $j=n$ and $L_{0,j,0}=a_{1}^{p}a_{2}^{jp}\in
x_{1}I^{n-1}$. Suppose $0<k\leq p$. Then
$L_{0,j,k}=a_{1}^{p}a_{2}^{jp+k(p-1)}=a_{1}^{p}(a_{2}^{p})^{j+k-1}a_{2}^{p-k}\in
x_1I^{n-1}$. Finally, if $k>p$, then
$L_{0,j,k}=a_{2}^{jp}(a_{1}a_{2}^{p-1})^{k}=
a_{1}^{p}[a_{1}^{k-p}a_{2}^{jp+k(p-1)}]=
a_{1}^{p}[a_{1}^{k-p}a_{2}^{(k-p)(p-1)}a_{2}^{jp+p(p-1)}]=
a_{1}^{p}(a_{2}^{p})^{j+p-1}(a_{1}a_{2}^{p-1})^{k-p}\in x_{1}I^{n-1}$.

Note that, if $p>2$, then $x_{2}^{*}$ is not a $\graded(I)$-sequence
because
\begin{eqnarray*}
(a_{1}^{p+2}a_{2}^{p-2})x_{2}=a_{1}^{p}(a_{1}a_{2}^{p-1})^{2}=x_{1}y^{2}\in
I^{3},
\end{eqnarray*}
where $(a_{1}^{p+2}a_{2}^{p-2})\in I\setminus I^{2}$.

Therefore, we can apply Theorem~\ref{main-theorem} and deduce that,
for all $n\geq 2$,
\begin{eqnarray*}
\left[ \frac{Q}{Q\langle n-1\rangle}\right]_{n}\cong
\frac{JI^{n-1}:y^{n}}{JI^{n-2}:y^{n-1}}.
\end{eqnarray*}
In particular, since $(JI^{p-2}:y^{p-1})\subsetneq
(JI^{p-1}:y^{p})=R$, then $\reltype(I)=\rednumber_J(I)+1$.

\noindent \underline{\sc Claim}: $(JI^{n-1}:y^{n})=(a_1^{p-n},a_2)$
for $2\leq n\leq p-1$. First note that, for all $2\leq n\leq p-1$,
$a_2\in (JI^{n-1}:y^{n})$ since
$a_2y^n=a_1^{n}a_2^{n(p-1)+1}=a_1a_2^p(a_1^{n-1}a_2^{(n-1)(p-1)})=
a_1a_2^p(a_1a_2^{p-1})^{n-1}\in JI^{n-1}$. Since $(JI^{n-1}:y^{n})$ is
generated by monomials on $a_1,a_2$, and $a_{2}\in (JI^{n-1}:y^{n})$,
there is only one possible remaining generator: a power of
$a_1$. Since $a_1^{p-n}y^n=a_1^{p-n} a_1^n a_2^{n(p-1)}=
a_1^pa_2^{n(p-1)}=a_2^{p-n}a_1^p(a_2^p)^{n-1}\in JI^{n-1}$, then
$a_{1}^{p-n}\in (JI^{n-1}:y^{n})$. However, and using again
\cite[Corollary~3]{ks}, one has $a_1^{p-n-1}\notin (JI^{n-1}:y^{n})$.

Hence $(Q/Q\langle p-1\rangle)_{p}\cong
(JI^{p-1}:y^{p})/(JI^{p-2}:y^{p-1})=R/(a_{1},a_{2})$ and, for $2\leq
n\leq p-1$,
\begin{eqnarray*}
\left[ \frac{Q}{Q\langle n-1\rangle}\right]_{n}\cong
\frac{JI^{n-1}:y^{n}}{JI^{n-2}:y^{n-1}}=
\frac{(a_{1}^{p-n},a_{2})}{(a_{1}^{p-n+1},a_{2})}\cong
\frac{(a_{1}^{p-n})}{(a_{1}^{p-n+1},a_{1}^{p-n}a_{2})}.
\end{eqnarray*}
In other words, for each $2\leq n\leq p$, $(Q/Q\langle
n-1\rangle)_{n}$ is generated by a single element that corresponds to
the class of $a_{1}^{p-n}$ ($1$ if $n=p$). To find this element,
consider the identity $a_1^{p-n}y^n=a_2^{p-n}a_1^p(a_2^p)^{n-1}$,
which induces the equation
$F_{n}(X_1,X_2,Y)=a_1^{p-n}Y^n-a_2^{p-n}X_1X_2^{n-1}\in Q_n$. Since
the isomorphism of Theorem \ref{main-theorem} sends the class of $F_n$
to the class of $F_n(0,0,1)=a_{1}^{p-n}$, we are done. By
Remark~\ref{min-gen-set-of-equations}, $F_n$ is in a minimal
generating set of $Q$, for $2\leq n\leq p$.

To finish, let us find the equations of degree one. Though this is
trivial, we sketch the proof here to show the similarity with the
greater degrees. As before, one shows that $J=(a_{1}^{p-1},a_2)$.
Using Remark~\ref{syzygies},
\begin{eqnarray*}
\koszul{1}{x_1,x_2,y}{R}\stackrel{\tilde{\sigma}_1}{\cong}
\frac{(J:y)}{J}=\frac{(a_1^{p-1},a_2)}{(a_1^p,a_2^p)}.
\end{eqnarray*}
Identify $Q_1$ with $\cycle{1}{x_1,x_2,y}{R}$ and $B_1=\langle
x_1Y-yX_1, x_2Y-yX_2,x_1X_2-x_2X_1\rangle$ with
$\boundary{1}{x_1,x_2,y}{R}$.  Then $Q_1/B_1$ is minimally generated
by the classes of the two equations corresponding to the classes of
$a_1^{p-1}$ and $a_2$. The identities $a_1^{p-1}y=a_2^{p-1}a_1^p$ and
$a_2y=a_1a_2^p$ induce the desired equations
$F_1(X_1,X_2,Y)=a_1^{p-1}Y-a_2^{p-1}X_1$ and
$G_1(X_1,X_2,Y)=a_2Y-a_1X_2\in Q_1$, since by $\tilde{\sigma}_1$ their
classes are sent to the classes of $F_1(0,0,1)=a_1^{p-1}$ and
$G_1(0,0,1)=a_2$. Clearly, $B_1\subset \langle F_1,G_1\rangle$ and
$F_1,G_1$ are a minimal generating set of $Q_1$.
\end{proof}

In the next example, although the initial forms $x_1^*,x_2^*$ are
zero-divisors in $\graded(I)$, where $I=(x_1,x_2,y)$, one may still
apply Theorem~\ref{main-theorem}, because the conditions of
Lemma~\ref{homologia-zero} ``$(J_{i-1}I^{n-1}:x_i)\cap
I^{n-1}/J_{i-1}I^{n-2}=0$, for all $i=1,\ldots ,s$, for all $n\geq
2$'', still hold (see Remark~\ref{weaker-hyp}).

\begin{example}\label{gen-classic} {\rm 
Let $(R,\fm)$ be a Noetherian local ring. Let $a_1,a_2$ an
$R$-sequence and let $p\geq 5$ be an odd integer. Set
$x_{1}=a_{1}^{p}$, $x_{2}=a_{2}^{p}$ and $y=a_{1}^2a_{2}^{p-2}$. Let
$I$ be the ideal generated by $x_{1},x_{2},y$ and set
$J=(x_1,x_2)$. Then $J$ is a reduction of $I$ with reduction number
$\rednumber_J(I)=p-1$ and $(x_1I^{n-1}:x_2)\cap I^{n-1}/x_1I^{n-2}=0$,
for all $n\geq 2$. However, the initial forms $x_1^*,x_2^*$ are
zero-divisors in $\graded(I)$. Applying Remark~\ref{weaker-hyp}, one
obtains the minimal generating set of the equations of $\rees(I)$:
\begin{itemize}
\item an equation $F_n(X_1,X_2,Y)=Y^p-X_1^2X_2^{p-2}$ of degree
  $p$;
\item an equation $F_n(X_1,X_2,Y)=a_2Y^n-a_1X_1X_2^{n-1}$ of degree
  $n=(p+1)/2$;
\item an equation
  $F_n(X_1,X_2,Y)=a_1^{p-2n}Y^n-a_2^{p-2n}X_1X_2^{n-1}$ of degree $n$,
  for each $n$, $2\leq n\leq (p-1)/2$;
\item two equations $F_1(X_1,X_2;Y)=a_1^{p-2}Y-a_2^{p-2}X_1$ and
  $G_1(X_1,X_2,Y)=a_2^2Y-a_1^2X_2$ of degree 1.
\end{itemize}
We leave the details to the reader.  }\end{example}

The next example shows that the weaker hypotheses in
Remark~\ref{different-version-theorem} are not enough to assure that
there is only one equation of $\rees(I)$ of maximum degree.

\begin{example}\label{not-enough} {\rm 
Let $k$ be a field and $R=k[X,Y]_{(X,Y)}/(XY,Y^{2})_{(X,Y)}$. Set $x$
and $y$ the classes of $X$ and $Y$ in $R$. Let $\fm=(x,y)$ be the
maximal ideal of $R$. Then the relation type of $\fm$ is 2 and there
are two equations of degree 2 in a minimal generating set of equations
of $\rees(\fm)$.  }\end{example}
\begin{proof}
Set $J=(x)$. Since $y^2=0\in J\fm$ and $y\not\in J$, then $\fm^2=J\fm$
and $J$ is a reduction of $\fm$ with reduction number $1$. Moreover,
since $(0:x)=(0:x^2)$, $x$ is a $d$-sequence. By
Remark~\ref{different-version-theorem}, $\reltype(\fm)\leq
\rednumber_J(\fm)+1=2$. Set $V=k[S,T]$ and let $\psi:V\to
\graded(\fm)$ be the presentation of $\graded(\fm)$ sending $S$ to
$x+\fm^2$ and $T$ to $y+\fm^2$. For $n\geq 2$, $\fm^n=(x^n)$.  Thus
$\fm^2/\fm^3$ is a $k$-vector space of dimension 1. Therefore
$\ker(\psi_2)\subset V_2$ must have dimension 2. In fact,
$\ker(\psi_2)=\langle ST,T^2\rangle$. Since $\ker(\psi_1)=0$, then
$E(\graded(\fm))_2$ is minimally generated by two elements. We finish
by using that $E(\fm)_2\cong E(\graded(\fm))_2$
(\cite[Proposition~3.3]{planas}).
\end{proof}

The next example shows that if the deviation $\mu(I)-\height(I)$ fails
to be equal to $1$, then there might be several equations of
$\rees(I)$ of maximum degree.

\begin{example}\label{exemple-m-h=2} {\rm 
Let $(R,\fm)$ be a two dimensional regular local ring, $p\geq 2$.
Then $\reltype(\fm^p)=2$ and there are $\binom{p}{2}$ equations of
degree 2 in a minimal generating set of equations of $\rees(\fm^p)$.
}\end{example}
\begin{proof}
Let $x,y$ be a regular system of parameters of $R$, $V=R[X,Y]$ and let
$\varphi:V\to\rees(\fm)$ be the presentation of $\rees(\fm)$ sending
$X$ to $xt$ and $Y$ to $yt$. Since $x,y$ is an $R$-sequence,
$\ker(\varphi)=\langle xY-yX\rangle$. Set $V(p)=R[X^p,X^{p-1}Y,\ldots
  ,XY^{p-1},Y^p]$ and $\varphi(p):V(p)\to \rees(\fm^p)$ the $p$-th
Veronese transform of $\varphi$. Note that
$\ker(\varphi(p)_n)=\ker(\varphi_{pn})=\langle xY-yX\rangle V_{pn-1}$.

Set $W=R[T_0,T_1,\ldots ,T_p]$ and let $\psi:W\to V(p)$ be the
polynomial presentation of $V(p)$ sending $T_i$ to $X^{p-i}Y^{i}$. It
is known (see e.g. \cite[Proposition~2.5]{jk}) that the kernel of
$\psi$ is the determinantal ideal generated by the $2\times 2$ minors
of the matrix $\cm$, where $\cm$ is
\begin{displaymath}
\left( \begin{array}{cccc} T_0 & T_1 & \ldots & T_{p-1} \\ T_1 & T_2 &
  \ldots & T_p \end{array} \right).
\end{displaymath}
In particular, $\ker(\psi)=\ker(\psi)\langle 2\rangle$.

Consider $\Phi=\varphi(p)\circ\psi:R[T_0,T_1,\ldots ,T_p]\to
\rees(\fm^p)$, the polynomial presentation of $\rees(\fm^p)$ sending
$T_i$ to $x^iy^{p-i}$ and let $Q=\ker(\Phi)$ be the ideal of equations
of $\rees(\fm^p)$. Let us see that $Q_n=\ker(\psi_n)+W_{n-1}Q_1$, for
all $n\geq 2$. Indeed, given $F\in Q_n$, since
$\ker(\varphi(p)_n)=\langle xY-yX\rangle V_{pn-1}$, one can find
$F_i\in W_{n-1}$ and $G_i\in Q_1$ such that $\psi _n(F)=\psi _n(\sum
F_iG_i)$. Therefore, $(F-\sum F_iG_i)\in \ker(\psi_n)$ and
$F\in\ker(\psi_n)+W_{n-1}Q_1$.

Since $\ker(\psi)=\ker(\psi)\langle 2\rangle$, then
$Q_n=\ker(\psi_n)+W_{n-1}Q_1=W_{n-2}\ker(\psi_2)+W_{n-1}Q_1\subseteq
W_{n-2}Q_2\subseteq Q_n$, for all $n\geq 2$. Therefore $Q=Q\langle
2\rangle$, $\reltype(\fm^n)=2$ and $\mu(Q_2/W_1Q_1)\leq\binom{p}{2}$.

Let us prove that $Q_2/W_1Q_1$ is minimally generated by
$\binom{p}{2}$ elements, which are precisely the classes of the
$2\times 2$ minors of $\cm$. Since $Q_1\subset \fm W_1$, then
$W_1Q_1\subset \fm W_2$. Setting $L=\ker(\psi_2)$, $M=W_2$ and
$N=V(p)_2$, we have
\begin{eqnarray*}
\frac{Q_2}{W_1Q_1}\otimes k\cong \frac{Q_2}{\fm
  Q_2+W_1Q_1}=\frac{\ker(\psi_2)+W_1Q_1} {\fm\ker(\psi_2)+W_1Q_1}
\cong \frac{L}{\fm L+(L\cap W_1Q_1)}.
\end{eqnarray*}
On the other hand, there is a natural epimorhism 
\begin{eqnarray*}
\frac{L}{\fm L+(L\cap W_1Q_1)}\to \frac{L}{L\cap \fm M},
\end{eqnarray*}
where $L/(L\cap\fm M)\cong (L+\fm M)/\fm M=\ker(\psi_2\otimes
1_k)$. Hence, 
\begin{eqnarray*}
\mu\left( \frac{Q_2}{W_1Q_1}\right) \geq \dim \frac{L}{(L\cap \fm
  M)}=\dim \ker(\psi_2\otimes 1_k),
\end{eqnarray*}
which clearly is $\binom{p}{2}$. Therefore $\mu((Q/Q\langle
1\rangle)_2)=\binom{p}{2}$ and, by
Remark~\ref{min-gen-set-of-equations}, there are $\binom{p}{2}$
equations of degree 2 in a minimal generating set of equations of
$\rees(\fm^p)$.

Note that $J=(x_1^p,x_2^p)$ is a reduction of $\fm^p$ with reduction
number $\rednumber_J(\fm^p)=1$, that $x_1^p,x_2^p$ is an $R$-sequence
and that $(x_1^p)^*$ is a $\graded(\fm^p)$-sequence.
\end{proof}

\bibliographystyle{amsplain}

\vspace{1cm}
\noindent \small Departament de Matemàtica Aplicada~1, Universitat Polit\`ecnica de Catalunya, Diagonal~647, ETSEIB, E-08028 Barcelona, Catalunya, Spain\\[0.3cm]
\small Email: \texttt{ferranmuinos@gmail.com}, \texttt{francesc.planas@upc.edu}.

\end{document}